\documentclass[11pt]{article}
\usepackage{amsmath}
\usepackage{amsfonts}
\usepackage{amssymb}
\usepackage{hyperref}
\usepackage{indentfirst}
\usepackage{mathrsfs}

\usepackage[top=1in,bottom=1in,left=1.25in,right=1.25in]{geometry}
\date{}
\allowdisplaybreaks

\begin{document}

\renewcommand{\baselinestretch}{1.2}
\renewcommand{\arraystretch}{1.0}

\title{\bf A Class of Quasitriangular Group-cograded Multiplier Hopf Algebras}
\author
{
  \textbf{Tao Yang} \footnote{Corresponding author. College of Science, Nanjing Agricultural University, Nanjing 210095, Jiangsu, China.
             E-mail: tao.yang@njau.edu.cn}, \,
  \textbf{Xuan Zhou} \footnote{Department of Mathematics, Jiangsu Second Normal University, Nanjing 210013, Jiangsu, China. E-mail: 20668964@qq.com},
  \,
  \textbf{Haixin Zhu} \footnote{College of Economics and Management, Nanjing Forestry University, Nanjing 210037, Jiangsu, China. E-mail: zhuhaixing@163.com}
}

\maketitle

\begin{center}
\begin{minipage}{12.cm}

 {\bf Abstract }
 For a multiplier Hopf algebra pairing $\langle A, B\rangle$, we construct a class of group-cograded multiplier Hopf algebras $D(A, B)$,
 generalizing the classical construction of finite dimensional Hopf algebras introduced by Panaite and Staic Mihai \cite{PS07}.
 Furthermore, if the multiplier Hopf algebra pairing admits a canonical multiplier in $M(B\otimes A)$
 we show the existence of quasitriangular structure on $D(A, B)$. As an application, some special cases and examples are provided.
\\

\\


\end{minipage}
\end{center}
\normalsize

\section{Introduction}
\def\theequation{\thesection.\arabic{equation}}
\setcounter{equation}{0}

 Recall from \cite{DVW05} that the motivating example for quasitriangular Hopf algebras is $H = U_{q}(\mathfrak{g})$,
 where $\mathfrak{g}$ is a finite-dimensional semisimple Lie algebra over the field $\mathbb{C}$ of complex numbers.
 In fact, $H$ is not quasitriangular in the strict sense of the definition,
 because the R-matrix lies in a completion of $H \otimes H$ rather than in $H \otimes H$ itself.
 The Hopf algebra $H = U_{q}(\mathfrak{g})$ is ''topologically'' quasitriangular.
 The explicit construction of the universal R-matrix is complicated.
 One approach with multiplier Hopf algebras gives a way to construct a generalized R-matrix in purely algebraic terms.
 The notion of a quasitriangular multiplier Hopf algebra is introduced in \cite{Zh99}.

 The concept of a group-cograded multiplier Hopf algebra was introduced by Abd El-hafez, Delvaux and Van Daele in \cite{ADV07}
 as a generalization of Hopf group-coalgebras introduced in \cite{T}.
 In \cite{DVW05}, the authors brought the results of quasitriangular Hopf group-coalgebras (as introduced by Turaev)
 to the more general framework of multiplier Hopf algebras, i.e., quasitriangular group-cograded multiplier Hopf algebras.

 In \cite{PS07}, the authors constructed a class of Hopf group-coalgebras by the so-called diagonal crossed product of a
 finite dimensional Hopf algebra $H$ and its duality $H^*$.
 Then one main question arises: Does the construction still holds for some infinite dimensional Hopf algebra?

 For the question, we first consider a more general case, and give a positive answer to this question.
 The main aim of this paper is to construct more examples
 of quasitriangular group-cograded multiplier Hopf algebras.

 The paper is organized in the following way.
 In section 2, we recall some notions which will be used in the following, such as
 multiplier Hopf algebras, quasitriangular group-cograded multiplier Hopf algebras and pairing.

 In section 3, let $A$ and $B$ be regular multiplier Hopf algebras with pairing $\langle A, B\rangle$.
 Then $D(A, B) = \bigoplus_{(\alpha, \beta)\in G} A \bowtie B_{(\alpha, \beta)}$
 is a $G$-cograded multiplier Hopf algebra, where $A \bowtie B_{(\alpha, \beta)}$ is the diagonal crossed product
 and $G = Aut_{Hopf}(B)\times Aut_{Hopf}(B)$ is a group with multiplication
 $(\alpha, \beta) \ast (\gamma, \delta) = (\alpha\gamma, \delta\gamma^{-1}\beta\gamma)$
 for $\alpha, \beta, \gamma, \delta \in Aut_{Hopf}(B)$.

 In section 4, we show in Theorem 4.3 that $D(A, B)$ constructed in the section 3 admits a quasitriangular structure if
 there is a canonical multiplier in $M(B\otimes A)$.

 In section 5, we also conclude by describing its applications and examples in the setting of Hopf algebras.

\section{Preliminaries}
\def\theequation{\thesection.\arabic{equation}}
\setcounter{equation}{0}

 We begin this section with a short introduction to multiplier Hopf algebras.

 Throughout this paper, all spaces we considered are over a fixed field $K$ (such as the field $\mathbb{C}$).
 Algebras may or may not have units, but should be always non-degenerate,
 i.e., the multiplication maps (viewed as bilinear forms) are non-degenerate.
 For an algebra $A$, the multiplier algebra $M(A)$ is defined as the largest algebra with unit
 in which $A$ is a dense ideal (see the appendix in \cite{V94}).

 Now, we recall the definition of a multiplier Hopf algebra (see \cite{V94} for details).
 A comultiplication on an algebra $A$ is a homomorphism $\Delta: A \longrightarrow M(A \otimes A)$
 such that $\Delta(a)(1 \otimes b)$ and $(a \otimes 1)\Delta(b)$ belong to $A\otimes A$ for all $a, b \in A$.
 We require $\Delta$ to be coassociative in the sense that
 \begin{eqnarray*}
 (a\otimes 1\otimes 1)(\Delta \otimes \iota)(\Delta(b)(1\otimes c))
 = (\iota \otimes \Delta)((a \otimes 1)\Delta(b))(1\otimes 1\otimes c)
 \end{eqnarray*}
 for all $a, b, c \in A$, where $\iota$ denotes the identity map.

 A pair $(A, \Delta)$ of a non-degenerate algebra $A$ with a comultiplication $\Delta$ is called a \emph{multiplier Hopf algebra},
 if the maps $T_{1}, T_{2}: A\otimes A \longrightarrow M(A\otimes A)$ defined by
 \begin{eqnarray}
 T_{1}(a\otimes b)=\Delta(a)(1 \otimes b), \qquad  T_{2}(a\otimes b)=(a \otimes 1)\Delta(b)
 \end{eqnarray}
 have range in $A\otimes A$ and are bijective.

 A multiplier Hopf algebra $(A, \Delta)$ is called \emph{regular} if $(A, \Delta^{cop})$ is also a multiplier Hopf algebra,
 where $\Delta^{cop}$ denotes the co-opposite comultiplication defined as $\Delta^{cop}=\tau \circ \Delta$ with $\tau$ the usual flip map
 from $A\otimes A$ to itself (and extended to $M(A\otimes A)$). In this case,
 $\Delta(a)(b \otimes 1) \mbox{ and } (1 \otimes a)\Delta(b) \in A \otimes A$
 for all $a, b\in A$.

 By Proposition 2.9 in \cite{V98},
 Multiplier Hopf algebra $(A, \Delta)$ is regular if and only if the antipode $S$ is bijective from $A$ to $A$.
 In this situation, the comultiplication is also determined by the bijective maps
 $T_{3}, T_{4}: A\otimes A \longrightarrow A\otimes A$ defined as follows
 \begin{eqnarray}
 && T_{3}(a \otimes b)=\Delta(a)(b \otimes 1), \qquad T_{4}(a\otimes b)=(1 \otimes a)\Delta(b)
 \end{eqnarray}
 for all $a, b\in A$.

 In this paper, all the multiplier hopf algebras we considered are regular.
 We will use the adapted Sweedler notation for regular multiplier Hopf algebras (see \cite{V08}).
 We will e.g., write $\sum a_{(1)} \otimes a_{(2)}b$ for $\Delta(a)(1 \otimes b)$
 and $\sum ab_{(1)} \otimes b_{(2)}$ for $(a \otimes 1)\Delta(b)$, sometimes we omit the $\sum$.

\subsection{Quasitriangular group-cograded multiplier Hopf algebras}

 The concept of a group-cograded multiplier Hopf algebra was introduced by Abd El-hafez, Delvaux and Van Daele in \cite{ADV07}
 as a generalization of Hopf group-coalgebras  introduced in \cite{T}.

 Let $(A, \Delta)$ be a multiplier Hopf algebra and $G$ a group.
 Assume that there is a family of (non-trivial) subalgebras $(A_{p})_{p\in G}$ of $A$ so that
 \begin{itemize}
 \item[(i)] $A=\bigoplus_{p\in G}A_{p}$ with $A_{p}A_{q}=0$ whenever $p,q\in G$ and $p\neq q$,
 \item[(ii)] $\Delta(A_{pq})(1\otimes A_{q})=A_{p}\otimes A_{q}$ and $(A_{p}\otimes 1)\Delta(A_{pq})=A_{p}\otimes  A_{q}$ for all $p, q \in G$.
 \end{itemize}
 Then $(A, \Delta)$ is called  a {\it $G$-cograded multiplier Hopf algebra}.
 The theory of group-cograded multiplier Hopf algebras was further developed.
 In particular in \cite{DVW05}, the authors study quasitriangular group-cograded multiplier Hopf algebras in the following sense:
 a $G$-cograded multiplier Hopf algebra with a crossing action $\xi$ is called {\it quasitriangular} if there is a multiplier
 $R=\sum _{\alpha, \beta \in G}R_{\alpha , \beta }$ with $R_{\alpha , \beta }
 \in M(A_{\alpha }\otimes A_{\beta })$ so that
 (1): $(\xi _p\otimes \xi _p)(R) = R$ for all $p\in G$, (2): $(\widetilde{\Delta}\otimes  \iota)(R)=R_{13}R_{23},
 (3): (\iota \otimes \Delta )(R)=R_{13}R_{12}$ and  (4):  $R\Delta (a)=(\widetilde{\Delta} )^{cop}(a)R$
 for all $p\in G$ and $a\in A$, where $\widetilde{\Delta}(a)(1\otimes a')=(\xi_{q^{-1}}\otimes \i)(\Delta(a)(1\otimes a'))$,
 for all $a\in A$ and $a'\in A_{q}$.

 \subsection{Multiplier Hopf algebra pairing}

 Start with two regular multiplier Hopf algebras $A$ and $B$ together with a non-degenerate bilinear map $\langle\cdot, \cdot\rangle$
 from $A\times B$ to $K$ satisfying certain properties. The main property is the comultiplication in $A$ is dual to the product in $B$ and vice versa.
 For more details, see \cite{DrV01}.

 For $a\in A$ and $b\in B$, we can define multipliers $a\blacktriangleright
 b$, $b\blacktriangleright a$, $a\blacktriangleleft b$ and $b\blacktriangleleft
 a$ in the following way. For $a'\in A$ and $b'\in B$, we have:
 $(b\blacktriangleright a)a'=\sum \langle a_{(2)}, b\rangle a_{(1)}a'$,
 \, $(a\blacktriangleright b)b'=\sum \langle a, b_{(2)} \rangle b_{(1)}b'$,
 \, $(a\blacktriangleleft b)a'=\sum \langle a_{(1)}, b \rangle a_{(2)}a'$
 and $(b\blacktriangleleft a)b'=\sum \langle a, b_{(1)} \rangle b_{(2)}b'$.
 The regularity conditions on the dual paring $\langle , \rangle$ say that the multipliers
 $b\blacktriangleright a$ and $a\blacktriangleleft b$ in $M(A)$
(resp. $a\blacktriangleright b$ and $b\blacktriangleleft a$ in $M(B)$)
 actually belong to $A$ (resp. $B$). For more details, see \cite{DV04}.

 We mention that $\langle S(a), b\rangle = \langle a, S(b)\rangle$, $\langle 1_{M(A)}, b\rangle = \varepsilon(b)$
 and $\langle a, 1_{M(B)}\rangle = \varepsilon(b)$.
 Sometimes without confusion we denote the unit $1_{M(A)}$ of $M(A)$ by $1$.
 We also use bilinear forms on the tensor products in the following way
 \begin{eqnarray*}
 \langle a\otimes a', b\otimes b'\rangle = \langle a, b\rangle \langle a',  b'\rangle, \quad
 \langle a\otimes b, b'\otimes a'\rangle = \langle a, b\rangle \langle a',  b'\rangle
 \end{eqnarray*}
 for all $a, a'\in A$ and $b, b'\in B$.
 These bilinear forms are non-degenerate and can be extended in a natural way to the multiplier algebra at one side.

\section{Diagonal crossed product of multiplier Hopf algebras}
\def\theequation{\thesection.\arabic{equation}}
\setcounter{equation}{0}

 Let $B$ be a multiplier Hopf algebra, we denote the group of multiplier Hopf automorphisms by $Aut_{Hopf}(B)$.
 Let $\alpha \in Aut_{Hopf}(B)$, by Lemma 3.3 in \cite{DV07} we have $\Delta\circ\alpha =(\alpha\otimes\alpha)\circ\Delta$,
 $\varepsilon\circ \alpha = \varepsilon$ and $S\circ\alpha = \alpha\circ S$.
 Denote $G = Aut_{Hopf}(B)\times Aut_{Hopf}(B)$, a group with multiplication
 \begin{eqnarray}
 (\alpha, \beta) \ast (\gamma, \delta) = (\alpha\gamma, \delta\gamma^{-1}\beta\gamma), \label{3.1}
 \end{eqnarray}
 The unit is $(\iota, \iota)$ and $(\alpha, \beta)^{-1} = (\alpha^{-1}, \alpha\beta^{-1}\alpha^{-1})$ (see \cite{YW11}).
 \\

 Firstly, we introduce the diagonal crossed product of regular multiplier Hopf algebras.
 On the other words, the construction of the diagonal crossed product in \cite{PS07} still holds for the multiplier Hopf algebras.
 However, we need to check that the diagonal crossed product is non-degenerate.
 \\

 \textbf{Definition \thesection.1}
 Let $A$ and $B$ be regular multiplier Hopf algebras with pairing $\langle A, B\rangle$.
 For $(\alpha, \beta) \in G$, we set $A \bowtie B_{(\alpha, \beta)} = A\otimes B$ as a vector space
 with a multiplication defined by the following formula:
 \begin{eqnarray}
 (a\bowtie b)(a'\bowtie b')
 = a \big( \alpha(b_{(1)}) \blacktriangleright a' \blacktriangleleft S^{-1}\beta(b_{(3)})  \big) \bowtie b_{(2)}b' \label{a}
 \end{eqnarray}
 for all $a, a'\in A$ and $b, b'\in B$.
 This multiplication is called the diagonal crossed product.
 \\

 \textbf{Remark}
 The diagonal crossed product (\ref{a}) is well-defined. Indeed, for $a'\in A$ there is an element $e\in B$ such that
 $e\blacktriangleright a' = a'$, therefore the right side of equation (\ref{a}) becomes
 $a \big( \alpha(b_{(1)} \alpha^{-1}(e)) \blacktriangleright a' \blacktriangleleft S^{-1}\beta(b_{(3)})  \big) \bowtie b_{(2)}b'$.
 $b_{(1)} \alpha^{-1}(e) \otimes b_{(2)}b' \otimes b_{(3)}
 = (\iota\otimes\Delta)\big(\Delta(b)(\alpha^{-1}(e)\otimes 1)\big)(1\otimes b'\otimes 1) \in B\otimes B\otimes B$,
 so (\ref{a}) is well-defined.
 \\

 \textbf{Proposition \thesection.2}
 Take the notations as above.
 Then $A \bowtie B_{(\alpha, \beta)}$ with the diagonal crossed product defined by (\ref{a}) is an associative and non-degenerate algebra.
 Moreover, the algebras $A$ and $B$ are subalgebras of $A \bowtie B_{(\alpha, \beta)}$ by the linear embedding
 $A\hookrightarrow A \bowtie B_{(\alpha, \beta)}$ and $B\hookrightarrow A \bowtie B_{(\alpha, \beta)}$ defined by
 $a\mapsto a\bowtie 1_{M(B)}$ and $b\mapsto 1_{M(A)} \bowtie b$, respectively.

 \emph{Proof}
 We define two linear maps $t_{1}, t_{2}: A \otimes B \longrightarrow A \otimes B$ by the formulas:
 $t_{1}(a\otimes b) = \alpha(b_{(1)}) \blacktriangleright a \otimes b_{(2)}$ and
 $t_{2}(a\otimes b) = a \blacktriangleleft \beta(b_{(2)}) \otimes b_{(1)}$.
 Then $t_{1}$ and $t_{2}$ are bijective with the inverse given by
 $t_{1}^{-1}(a\otimes b) = S^{-1}\alpha(b_{(1)}) \blacktriangleright a \otimes b_{(2)}$ and
 $t_{2}^{-1}(a\otimes b) = a \blacktriangleleft S^{-1}\beta(b_{(2)}) \otimes b_{(1)}$ respectively.

 Let $T=t_{1}\circ t_{2}^{-1} \circ \tau$, then we have
 \begin{eqnarray*}
 T(b\otimes a') = \big( \alpha(b_{(1)}) \blacktriangleright a' \blacktriangleleft S^{-1}\beta(b_{(3)})  \big) \bowtie b_{(2)}
 \end{eqnarray*}
 is bijective. In this case the diagonal crossed product becomes the twisted tensor product in the sense of Delvaux \cite{D03},
 i.e., $(a\bowtie b)(a'\bowtie b') = (m_{A}\otimes m_{B})(\iota\otimes T\otimes \iota)(a\otimes b\otimes a'\otimes b')$.
 Then by Proposition 1.1 in \cite{D03}, the diagonal crossed product on $A \bowtie B_{(\alpha, \beta)}$ is non-degenerate.

 For the associativity and the rest of this proposition, it is straightforward.
 $\hfill \blacksquare$
 \\

 \textbf{Remark}
 (1) The product of $A \bowtie B_{(\alpha, \beta)}$ is non-degenerate, so we can get the multiplier Hopf algebra $M(A \bowtie B_{(\alpha, \beta)})$
 and obviously $1_{M(A)}\bowtie 1_{M(B)}$ is its unit.

 (2) By the 'cover technique' introduced in \cite{V08}, the product of $A \bowtie B_{(\alpha, \beta)}$
 can be written in adapted Sweedler notation:
 \begin{eqnarray*}
 (a\bowtie b)(a'\bowtie b') = \langle a'_{(1)}, S^{-1}\beta(b_{(3)})\rangle  \langle a'_{(3)}, \alpha(b_{(1)})\rangle ( aa'_{(2)}\bowtie b_{(2)}b').
 \end{eqnarray*}

 In particular, If $B$ is finite dimensional, then $B$ is a Hopf algebra. Let $A=H^*$,
 then the formula (\ref{a}) is just the diagonal crossed product introduced in \cite{PS07}.

 (3) As in Section 2.3 in \cite{DV04} the commutation rule in $A \bowtie B_{(\alpha, \beta)}$ canbe written as
 \begin{eqnarray}
 \langle a_{(1)}, b_{(2)}\rangle \big(1\bowtie \beta^{-1}(b_{(1)}) \big)(a_{(2)}x \bowtie y)
 = \langle  a_{(2)}, \alpha\beta^{-1}(b_{(1)})\rangle \big( a_{(1)}\bowtie \beta^{-1}(b_{(2)}) \big)  (x \bowtie y) \label{b}
 \end{eqnarray}
 for $a\in A$, $b\in B$ and $x\bowtie y\in A \bowtie B_{(\alpha, \beta)}$.
 \\

 In what follows, let $D(A, B) = \bigoplus_{(\alpha, \beta)\in G} A \bowtie B_{(\alpha, \beta)}$.
 Then we have the main results of this section: there exists a multiplier Hopf algebra structure on $D(A, B)$,
 which generalizes the classical construction of finite dimensional Hopf algebras by Panaite and Staic Mihai in \cite{PS07}.
 This construction is different from what introduced in \cite{YW11a}.
 \\

 \textbf{Theorem \thesection.3}
 Let $A$ and $B$ be regular multiplier Hopf algebras with pairing $\langle A, B\rangle$.
 Then $D(A, B) = \bigoplus_{(\alpha, \beta)\in G} A \bowtie B_{(\alpha, \beta)}$
 is a $G$-cograded multiplier Hopf algebra with the following strucrures:
 \begin{itemize}
 \item For any $(\alpha, \beta)\in G$,  the multiplication of $A \bowtie B_{(\alpha, \beta)}$ is given by Definition \thesection.1.

 \item The comultiplication on $D(A, B)$ is given by:
 \begin{eqnarray*}
 && \Delta_{(\alpha, \beta), (\gamma, \delta)}:
 A \bowtie B_{(\alpha, \beta)\ast (\gamma, \delta)} \longrightarrow A \bowtie B_{(\alpha, \beta)} \otimes A \bowtie B_{(\gamma, \delta)}, \\
 && \Delta_{(\alpha, \beta), (\gamma, \delta)} (a\bowtie b) = \Delta^{cop}(a)(\gamma\otimes \gamma^{-1}\beta\gamma)\Delta(b).
 \end{eqnarray*}

 \item The counit $\varepsilon_{D(A, B)}$ on $A \bowtie B_{(\iota, \iota)}$ is given by
 $\varepsilon_{D(A, B)}(a\bowtie b) = \varepsilon_{A}(a)\varepsilon_{B}(b)$.

 \item For any $(\alpha, \beta)\in G$, the antipode is given by
 \begin{eqnarray*}
 && S: A \bowtie B_{(\alpha, \beta)} \longrightarrow A \bowtie B_{(\alpha, \beta)^{-1}}, \\
 && S_{(\alpha, \beta)}(a\bowtie b) = T(\alpha\beta S(b) \otimes S^{-1}(a))
    \mbox{ in } A \bowtie B_{(\alpha, \beta)^{-1}} = A \bowtie B_{(\alpha^{-1}, \alpha\beta^{-1}\alpha^{-1})}.
 \end{eqnarray*}
 \end{itemize}

 \emph{Proof}
 It is easy to check that $\varepsilon_{D(A, B)}$ is a counit of $D(A, B)$.
 Similar to the Drinfel'd double for group-cograded multiplier Hopf algebras introduced in \cite{DV07},
 $\Delta_{(\alpha, \beta), (\gamma, \delta)} (a\bowtie b) \big(1_{D(A, B)} \otimes (a'\bowtie b') \big)
 \in A \bowtie B_{(\alpha, \beta)} \otimes A \bowtie B_{(\gamma, \delta)}$ and
 $\big((a''\bowtie b'') \otimes 1_{D(A, B)}\big) \Delta_{(\alpha, \beta), (\gamma, \delta)} (a\bowtie b)
 \in A \bowtie B_{(\alpha, \beta)} \otimes A \bowtie B_{(\gamma, \delta)}$
 for any $a\bowtie b \in A \bowtie B_{(\alpha, \beta)\ast (\gamma, \delta)}$,
 $a'\bowtie b' \in A \bowtie B_{(\gamma, \delta)}$ and $a''\bowtie b'' \in A \bowtie B_{(\alpha, \beta)}$.

 For the coassociativity, it is straightforward. Next, let us check that $\Delta_{(\alpha, \beta), (\gamma, \delta)}$ is multiplicative,
 i.e., $\Delta_{(\alpha, \beta), (\gamma, \delta)} \big((a\bowtie b)(a'\bowtie b') \big)
 = \Delta_{(\alpha, \beta), (\gamma, \delta)} (a\bowtie b) \Delta_{(\alpha, \beta), (\gamma, \delta)} (a'\bowtie b')$.
 Indeed, for any $a, a', a''\in A$ and $b, b', b''\in B$,
 \begin{eqnarray*}
  && \big((a''\bowtie 1_{M(B)})\otimes 1_{D(A, B)} \big) \Delta_{(\alpha, \beta), (\gamma, \delta)} \big((a\bowtie b)(a'\bowtie b') \big)
     \big(1_{D(A, B)}\otimes (1_{M(A)}\bowtie b'') \big) \\
  &=&  \langle a'_{(1)}, S^{-1}\delta\gamma^{-1}\beta\gamma (b_{(3)})\rangle  \langle a'_{(3)}, \alpha\gamma(b_{(1)})\rangle \\
    && \big( (a''\bowtie 1_{M(B)})\otimes 1_{D(A, B)} \big)
   \Delta_{(\alpha, \beta), (\gamma, \delta)} \big(aa'_{(2)}\bowtie b_{(2)}b' \big)
     \big(1_{D(A, B)}\otimes (1_{M(A)}\bowtie b'') \big)\\
  &=&  \langle a'_{(1)}, S^{-1}\delta\gamma^{-1}\beta\gamma (b_{(4)})\rangle  \langle a'_{(4)}, \alpha\gamma(b_{(1)})\rangle \\
     && \big(a''a_{(2)}a'_{(3)}\bowtie \gamma(b_{(2)}b'_{(1)}) \otimes a_{(1)}a'_{(2)}\bowtie \gamma^{-1}\beta\gamma(b_{(3)}b'_{(2)})b'' \big),
 \end{eqnarray*}
 and
 \begin{eqnarray*}
  && \big( (a''\bowtie 1_{M(B)})\otimes 1_{D(A, B)} \big) \Delta_{(\alpha, \beta), (\gamma, \delta)} (a\bowtie b)
      \Delta_{(\alpha, \beta), (\gamma, \delta)} (a' \bowtie b') \\
     && \big(1_{D(A, B)}\otimes (1_{M(A)}\bowtie b'') \big)  \\
  &=&  (a''a_{(2)}\bowtie \gamma(b_{(1)})) (a'_{(2)}\bowtie \gamma(b'_{(1)}) )
     \otimes (a_{(1)}\bowtie \gamma^{-1}\beta\gamma(b_{(2)})) (a'_{(1)}\bowtie \gamma^{-1}\beta\gamma(b'_{(2)})b'' )  \\
  &=&  \underline{\langle a'_{(4)}, S^{-1}\beta\gamma (b_{(3)})\rangle}  \langle a'_{(6)}, \alpha\gamma(b_{(1)})\rangle
      (a''a_{(2)}a'_{(5)} \bowtie \gamma(b_{(2)})b'_{(1)}) \\
      && \otimes  \langle a'_{(1)}, S^{-1}\delta\gamma^{-1}\beta\gamma (h_{(6)})\rangle
      \underline{\langle a'_{(3)}, \gamma\gamma^{-1}\beta\gamma(b_{(4)})\rangle}
      (a_{(1)}a'_{(2)} \bowtie \gamma^{-1}\beta\gamma(b_{(5)})b'_{(2)} b'' ) \\
  &=& \langle a'_{(1)}, S^{-1}\delta\gamma^{-1}\beta\gamma (b_{(4)})\rangle  \langle a'_{(4)}, \alpha\gamma(b_{(1)})\rangle \\
     && \big(a''a_{(2)}a'_{(3)}\bowtie \gamma(b_{(2)}b'_{(1)}) \otimes a_{(1)}a'_{(2)}\bowtie \gamma^{-1}\beta\gamma(b_{(3)}b'_{(2)})b'' \big).
 \end{eqnarray*}

 Because the $T$ is bijective, it is easy to get that the antipode $S$ is bijective. Also we have
 \begin{eqnarray*}
  S_{(\alpha, \beta)}(a\bowtie b)
  &=& \big(\beta S(b_{(3)}) \blacktriangleright S^{-1}(a) \blacktriangleleft \alpha(b_{(1)})\big) \bowtie \alpha\beta S(b_{(2)}) \\
  &=& \langle S^{-1}(a_{(3)}), \alpha(b_{(1)}) \rangle\langle S^{-1}(a_{(1)}), \beta S(b_{(3)}) \rangle
      \big( S^{-1}(a_{(2)})\bowtie \alpha\beta S(b_{(2)}) \big).
 \end{eqnarray*}
 It is straightforward to check that $S$ defined above is an algebra anti-isomorphism, i.e.,
 $S_{(\alpha, \beta)} \big((a\bowtie b) (a'\bowtie b')\big) = S_{(\alpha, \beta)}  (a'\bowtie b') S_{(\alpha, \beta)}  (a\bowtie b)$.
 In fact,
 \begin{eqnarray*}
  && S_{(\alpha, \beta)} \big((a\bowtie b) (a'\bowtie b')\big) \\
  &=& \langle a'_{(1)}, S^{-1}\beta(b_{(3)})\rangle  \langle a'_{(3)}, \alpha(b_{(1)})\rangle \big(S_{(\alpha, \beta)}
       (aa'_{(2)}\bowtie b_{(2)}b') \big) \\
  &=& \langle a'_{(1)}, S^{-1}\beta(b_{(3)})\rangle  \langle a'_{(3)}, \alpha(b_{(1)})\rangle
       \langle S^{-1}(a_{(3)}a'_{(4)}), \alpha(b_{(2)}b'_{(1)})\rangle \\
       && \langle S^{-1}(a_{(1)}a'_{(2)}), \beta S(b_{(4)}b'_{(3)})\rangle
       \big( S^{-1}(a_{(2)}a'_{(3)})\bowtie \alpha\beta S(b_{(3)}b'_{(2)})\big) \\
  &=& \langle a'_{(1)}, S^{-1}\beta(b_{(7)})\rangle  \langle a'_{(5)}, \alpha(b_{(1)})\rangle
       \langle S^{-1}(a'_{(4)}), \alpha(b_{(2)}b'_{(1)})\rangle  \langle S^{-1}(a_{(3)}), \alpha(b_{(3)}b'_{(2)})\rangle \\
       && \langle S^{-1}(a'_{(2)}), \beta S(b_{(6)}b'_{(5)})\rangle \langle S^{-1}(a_{(1)}), \beta S(b_{(5)}b'_{(4)})\rangle
       \big(S^{-1}(a_{(2)}a'_{(3)})\bowtie \alpha\beta S(b_{(4)}b'_{(3)})\big) \\
  &=& \langle a'_{(1)}, \beta(b'_{(5)})\rangle \langle a'_{(3)}, \alpha S^{-1}(b'_{(1)})\rangle
      \langle S^{-1}(a_{(3)}), \alpha(b_{(1)}b'_{(2)})\rangle \langle a_{(1)}, \beta (b_{(3)}b'_{(4)})\rangle \\
       && \big(S^{-1}(a_{(2)}a'_{(2)})\bowtie \alpha\beta S(b_{(2)}b'_{(3)})\big),
 \end{eqnarray*}
 and
 \begin{eqnarray*}
  && S_{(\alpha, \beta)} (a'\bowtie b') S_{(\alpha, \beta)} (a\bowtie b) \\
  &=& \langle S^{-1}(a'_{(3)}), \alpha(b'_{(1)}) \rangle\langle S^{-1}(a'_{(1)}), \beta S(b'_{(3)}) \rangle
      \big( S^{-1}(a'_{(2)})\bowtie \alpha\beta S(b'_{(2)}) \big) \\
    && \langle S^{-1}(a_{(3)}), \alpha(b_{(1)}) \rangle\langle S^{-1}(a_{(1)}), \beta S(b_{(3)}) \rangle
      \big( S^{-1}(a_{(2)})\bowtie \alpha\beta S(b_{(2)}) \big) \\
  &=& \langle S^{-1}(a'_{(3)}), \alpha(b'_{(1)}) \rangle\langle S^{-1}(a'_{(1)}), \beta S(b'_{(5)}) \rangle
      \langle S^{-1}(a_{(3)}), \alpha(b_{(1)}) \rangle\langle S^{-1}(a_{(1)}), \beta S(b_{(3)}) \rangle  \\
    && \langle S^{-1}(a_{(4)}), \alpha(b'_{(2)}) \rangle\langle S^{-1}(a_{(2)}), \beta S(b'_{(4)}) \rangle
       \big(S^{-1}(a_{(3)}a'_{(2)})\bowtie \alpha\beta S(b_{(2)}b'_{(3)})\big)\\
  &=& \langle a'_{(1)}, \beta(b'_{(5)})\rangle \langle a'_{(3)}, \alpha S^{-1}(b'_{(1)})\rangle
      \langle S^{-1}(a_{(3)}), \alpha(b_{(1)}b'_{(2)})\rangle \langle a_{(1)}, \beta (b_{(3)}b'_{(4)})\rangle \\
       && \big(S^{-1}(a_{(2)}a'_{(2)})\bowtie \alpha\beta S(b_{(2)}b'_{(3)})\big).
 \end{eqnarray*}

 Finally, we want to verify the following axiom: for $a\bowtie b \in A \bowtie B_{(\iota, \iota)}$,
 and $a'\bowtie b' \in A \bowtie B_{(\alpha, \beta)}$,
 \begin{eqnarray*}
  && m_{(\alpha, \beta)} \big(S_{(\alpha, \beta)^{-1}} \otimes \iota_{(\alpha, \beta)} \big)
  \big(\Delta_{(\alpha, \beta)^{-1}, (\alpha, \beta)} (a\bowtie b) (1 \otimes a'\bowtie b') \big)
  = \varepsilon_{D(A, B)}(a\bowtie b)(a'\bowtie b'),\\
  && m_{(\alpha, \beta)} \big(\iota_{(\alpha, \beta)} \otimes S_{(\alpha, \beta)^{-1}} \big)
  \big((a'\bowtie b' \otimes 1) \Delta_{(\alpha, \beta), (\alpha, \beta)^{-1}} (a\bowtie b) \big)
  = \varepsilon_{D(A, B)}(a\bowtie b)(a'\bowtie b').
 \end{eqnarray*}
 Here we only check the first equation, the second one is similar.
 \begin{eqnarray*}
  && m_{(\alpha, \beta)} \big(S_{(\alpha, \beta)^{-1}} \otimes \iota_{(\alpha, \beta)} \big)
    \big(\Delta_{(\alpha, \beta)^{-1}, (\alpha, \beta)} (a\bowtie b) (1 \otimes a'\bowtie b') \big) \\
 &=&   m_{(\alpha, \beta)} \big(S_{(\alpha, \beta)^{-1}} \otimes \iota_{(\alpha, \beta)} \big)
    \big(\Delta^{cop}(a)(\alpha\otimes\beta^{-1})\Delta(b) (1_{(\alpha, \beta)^{-1}} \otimes a'\bowtie b') \big) \\
 &=& S_{(\alpha^{-1}, \alpha\beta^{-1}\alpha^{-1})} \big( a_{(2)}\bowtie \alpha(b_{(1)}) \big)
    \big( a_{(1)}(\alpha\beta^{-1}(b_{(2)})\blacktriangleright a' \blacktriangleleft S^{-1}(b_{(4)})) \bowtie \beta^{-1}(b_{(3)})b' \big) \\
 &=& [\alpha\beta^{-1}S(b_{(3)})\blacktriangleright S^{-1}(a_{(2)}) \blacktriangleleft b_{(1)} \bowtie \beta^{-1}S(b_{(2)})] \\
   && [p_{(1)}\big(\alpha\beta^{-1}(h_{(4)})\blacktriangleright q \blacktriangleleft S^{-1}(h_{(6)})\big) \bowtie \beta^{-1}(h_{(5)})l] \\
 &=& [\alpha\beta^{-1}S(b_{(5)})\blacktriangleright S^{-1}(a_{(2)}) \blacktriangleleft b_{(1)}] \\
   &&  \Big[ \alpha\beta^{-1}S(b_{(4)})\blacktriangleright
     \Big( a_{(1)}\big(\alpha\beta^{-1}(b_{(6)})\blacktriangleright a' \blacktriangleleft S^{-1}(b_{(8)})\big)\Big)
     \blacktriangleleft b_{(2)} \Big] \bowtie \beta^{-1}S(b_{(3)})\beta^{-1}(b_{(7)})b' \\
 &=& [\alpha\beta^{-1}S(b_{(7)})\blacktriangleright S^{-1}(a_{(2)}) \blacktriangleleft b_{(1)}]
     \Big[ \alpha\beta^{-1}S(b_{(6)})\blacktriangleright a_{(1)}  \blacktriangleleft b_{(2)} \Big]\\
   &&  \Big[ \alpha\beta^{-1}S(b_{(5)})\blacktriangleright
     \Big(\alpha\beta^{-1}(b_{(8)})\blacktriangleright a' \blacktriangleleft S^{-1}(b_{(10)})\Big)
     \blacktriangleleft b_{(3)} \Big] \bowtie \beta^{-1}S(b_{(4)})\beta^{-1}(b_{(9)})b' \\
 &=& \Big[ \alpha\beta^{-1}S(b_{(5)})\blacktriangleright S^{-1}(a_{(2)})a_{(1)}  \blacktriangleleft b_{(1)} \Big]\\
   &&  \Big[ \alpha\beta^{-1}S(b_{(4)})\blacktriangleright
     \Big(\alpha\beta^{-1}(b_{(6)})\blacktriangleright a' \blacktriangleleft S^{-1}(b_{(8)})\Big)
     \blacktriangleleft b_{(2)} \Big] \bowtie \beta^{-1}S(b_{(3)})\beta^{-1}(b_{(7)})b' \\
 &=& \varepsilon(a)\varepsilon(b) a'\bowtie b'
 = \varepsilon_{D(A, B)}(a\bowtie b)(a'\bowtie b').
 \end{eqnarray*}
 Therefore, by Theorem 2.5 in \cite{ADV07} $D(A, B) = \bigoplus_{(\alpha, \beta)\in G} A \bowtie B_{(\alpha, \beta)}$
 is a regular $G$-cograded multiplier Hopf algebra.
 $\hfill \blacksquare$
 \\

 \textbf{Remark} Let $\pi$ be a subgroup of $Aut(B)$, then we also can construct the group $G' = \pi \times \pi$ by
 the product (\ref{3.1}). Then we can similarly obtain a group-cograded multiplier Hopf algebra over $G'$.
 \\

 \textbf{Example \thesection.4} Let $H$ be an infinite group. The Drinfel'd double $D(H)$ is a multiplier Hopf algebra
 rather than a usual Hopf algebra. Set $B=D(H)$, $A = \widehat{D(H)}$.
 Then the multiplier Hopf algebra structure on $B$ is given by
 \begin{eqnarray*}
 && (\delta_{p}\propto h)(\delta_{q}\propto l) = \delta_{p} \delta_{hqh^{-1}} \propto hl, \\
 && \Delta(\delta_{p}\propto h) = \sum_{s\in H} (\delta_{s^{-1}p}\propto h) \otimes (\delta_{s}\propto h), \\
 && \varepsilon (\delta_{p}\propto h) = \delta_{p, e}, \\
 && S(\delta_{p}\propto h) = \delta_{h^{-1}ph}\propto h^{-1},
 \end{eqnarray*}
 and the multiplier Hopf algebra structure on $A$ is given by
 \begin{eqnarray*}
 && (h\propto\delta_{p})(l\propto\delta_{q}) = lh \propto \delta_{p} \delta_{q}, \\
 && \Delta(h\propto\delta_{p}) = \sum_{t\in H} (h\propto\delta_{t}) \otimes (t^{-1}ht\propto\delta_{t^{-1}p}), \\
 && \varepsilon (h\propto\delta_{p}) = \delta_{p, e}, \\
 && S(h\propto\delta_{p}) = p^{-1}h^{-1}p\propto\delta_{p^{-1}}.
 \end{eqnarray*}
 Let $\alpha\in H$, define $\alpha(\delta_{p}\propto h) = \delta_{\alpha p\alpha^{-1}}\propto \alpha h\alpha^{-1} $,
 then $\alpha\in Aut_{Hopf}(D(H))$. By Theorem \thesection.3,
 $\mathcal{D}(D(H)) = \bigoplus_{(\alpha, \beta)\in G} \widehat{D(H)} \bowtie D(H)_{(\alpha, \beta)}$
 is a $G$-cograded multiplier Hopf algebra with the following strucrures:
 \begin{itemize}
 \item For any $(\alpha, \beta)\in G$,  the multiplication of $\widehat{D(H)} \bowtie D(H)_{(\alpha, \beta)}$ is given by
 \begin{eqnarray*}
 \Big( (1\bowtie(\delta_{p}\propto h) \Big)\Big( (l\propto \delta_{q})\bowtie 1 \Big)
 &=& (\beta h \beta^{-1} l \beta h^{-1}\beta^{-1}\propto \delta_{\beta h \beta^{-1} q \alpha h^{-1} \alpha^{-1} }) \\
  && \bowtie (\delta_{\alpha^{-1} l^{-1} \alpha ph\beta^{-1} l h^{-1}}\propto h).
 \end{eqnarray*}

 \item The comultiplication $\Delta_{(\alpha, \beta), (\gamma, \delta)}:
 \widehat{D(H)} \bowtie D(H)_{(\alpha, \beta)\ast (\gamma, \delta)}
    \longrightarrow \widehat{D(H)} \bowtie D(H)_{(\alpha, \beta)} \otimes \widehat{D(H)} \bowtie D(H)_{(\gamma, \delta)}$ is given by:
 \begin{eqnarray*}
 && \Delta_{(\alpha, \beta), (\gamma, \delta)} \big( (l\propto \delta_{q}) \bowtie (\delta_{p}\propto h) \big) \\
 &=& \sum_{s, t\in H} (t^{-1}l t \propto \delta_{t^{-1}q})\bowtie (\delta_{\gamma s \gamma^{-1}}\propto h) \\
 && \otimes
  (l\propto \delta_{t})\bowtie (\delta_{\gamma^{-1}\beta\gamma s^{-1}p \gamma^{-1}\beta^{-1}\gamma}\propto \gamma^{-1}\beta\gamma h
    \gamma^{-1}\beta^{-1}\gamma).
 \end{eqnarray*}

 \item The counit $\varepsilon_{\mathcal{D}(D(H))}$ on $\widehat{D(H)} \bowtie D(H)_{(\iota, \iota)}$ is given by
 \begin{eqnarray*}
 \varepsilon_{D(A, B)}\big( (l\propto \delta_{q}) \bowtie (\delta_{p}\propto h) \big) = \delta_{p, e}\delta_{q, e}.
 \end{eqnarray*}

 \item For any $(\alpha, \beta)\in G$, the antipode is given by $S: \widehat{D(H)} \bowtie D(H)_{(\alpha, \beta)} \longrightarrow \widehat{D(H)} \bowtie D(H)_{(\alpha, \beta)^{-1}}$
 \begin{eqnarray*}
 && S_{(\alpha, \beta)} \big( (l\propto \delta_{q}) \bowtie (\delta_{p}\propto h) \big) \\
 &=& (\alpha h^{-1} \alpha^{-1} q^{-1} lq\alpha h\alpha^{-1} \propto \delta_{\alpha h^{-1}\alpha^{-1} q^{-1} \beta h \beta^{-1}}) \\
 && \bowtie (\delta_{\alpha q^{-1}l^{-1}q\beta h^{-1} p \alpha^{-1} q^{-1} lq\alpha\beta h\beta^{-1}\alpha^{-1}}
    \propto \alpha\beta h^{-1}\beta^{-1}\alpha^{-1}).
 \end{eqnarray*}
 \end{itemize}

\section{Quasitriangular structures}
\def\theequation{\thesection.\arabic{equation}}
\setcounter{equation}{0}

 To construct quasitriangular structure on the $G$-cograded multiplier Hopf algebra established as before,
 we first study crossing actions as follows.

 \textbf{Proposition \thesection.1}
 With the notations as before. Then a crossing action $\xi: G \longrightarrow Aut(D(A, B))$ is given by
 \begin{eqnarray*}
  && \xi_{(\alpha, \beta)}^{(\gamma, \delta)}: A \bowtie B_{(\gamma, \delta)}
  \longrightarrow A \bowtie B_{(\alpha, \beta)\ast(\gamma, \delta)\ast(\alpha, \beta)^{-1}}
  = A \bowtie B_{(\alpha\gamma\alpha^{-1}, \alpha\beta^{-1}\delta\gamma^{-1}\beta\gamma\alpha^{-1})}, \\
  && \xi_{(\alpha, \beta)}^{(\gamma, \delta)}(a\bowtie b) = a\circ\beta\alpha^{-1} \bowtie \alpha\gamma^{-1}\beta^{-1}\gamma(b).
 \end{eqnarray*}

 \emph{Proof}
 First $\xi_{(\alpha, \beta)}^{(\gamma, \delta)}$ is an algebra morphism. Indeed,
 \begin{eqnarray*}
 && \xi_{(\alpha, \beta)}^{(\gamma, \delta)}(a\bowtie b) \xi_{(\alpha, \beta)}^{(\gamma, \delta)} (a'\bowtie b') \\
 &=& \big(a\circ\beta\alpha^{-1} \bowtie \alpha\gamma^{-1}\beta^{-1}\gamma(b) \big)
     \big(a'\circ\beta\alpha^{-1} \bowtie \alpha\gamma^{-1}\beta^{-1}\gamma(b') \big) \\
 &=& \langle a'_{(1)}\circ\beta\alpha^{-1}, S^{-1}\alpha\beta^{-1}\delta\gamma^{-1}\beta\gamma\alpha^{-1} \cdot
     \alpha\gamma^{-1}\beta^{-1}\gamma(b_{(3)})\rangle \\
   &&  \langle a'_{(3)}\circ\beta\alpha^{-1}, \alpha\gamma\alpha^{-1}\cdot \alpha\gamma^{-1}\beta^{-1}\gamma (b_{(1)})\rangle
     \big( (aa'_{(2)})\circ\beta\alpha^{-1} \bowtie \alpha\gamma^{-1}\beta^{-1}\gamma(b_{(2)}b') \big)\\
 &=& \langle a'_{(1)}, S^{-1}\delta (b_{(3)})\rangle
     \langle a'_{(3)}, \gamma (b_{(1)})\rangle
     \big( (aa_{(2)})\circ\beta\alpha^{-1} \bowtie \alpha\gamma^{-1}\beta^{-1}\gamma(b_{(2)}b') \big)\\
 &=& \xi_{(\alpha, \beta)}^{(\gamma, \delta)}\Big( \langle a'_{(1)}, S^{-1}\delta (b_{(3)})\rangle
     \langle a'_{(3)}, \gamma (b_{(1)})\rangle aa'_{(2)}\bowtie b_{(2)}b' \Big) \\
 &=& \xi_{(\alpha, \beta)}^{(\gamma, \delta)}\big((a\bowtie b)(a'\bowtie b')\big).
 \end{eqnarray*}
 Moreover $\alpha, \beta, \gamma, \delta \in Aut_{Hopf}(B)$ are bijective, then $\xi_{(\alpha, \beta)}^{(\gamma, \delta)}$ is an algebra isomorphism.

 Then it is straightforward to check that $\xi$ respects the comultiplication, i.e.,
 \begin{eqnarray*}
 && \Delta_{(\alpha, \beta)\ast(\gamma, \delta)\ast(\alpha, \beta)^{-1}, (\alpha, \beta)\ast(\mu, \nu)\ast (\alpha, \beta)^{-1}} \circ \xi_{(\alpha, \beta)}^{(\gamma, \delta)\ast (\mu, \nu)}
 = (\xi_{(\alpha, \beta)}^{(\gamma, \delta)} \otimes \xi_{(\alpha, \beta)}^{(\mu, \nu)})\circ \Delta_{(\gamma, \delta), (\mu, \nu)}.
 \end{eqnarray*}
 Indeed, for $a\bowtie b\in \mathcal{A}_{(\gamma, \delta)\ast (\mu, \nu)}$,
 \begin{eqnarray*}
 && \Delta_{(\alpha, \beta)\ast(\gamma, \delta)\ast(\alpha, \beta)^{-1}, (\alpha, \beta)\ast(\mu, \nu)\ast (\alpha, \beta)^{-1}} \circ \xi_{(\alpha, \beta)}^{(\gamma, \delta)\ast (\mu, \nu)} (a\bowtie b) \\
 &=& \Delta_{(\alpha\gamma\alpha^{-1}, \alpha\beta^{-1}\delta\gamma^{-1}\beta\gamma\alpha^{-1}),
    (\alpha\mu\alpha^{-1}, \alpha\beta^{-1}\nu\mu^{-1}\beta\mu\alpha^{-1})}
    \xi_{(\alpha, \beta)}^{(\gamma\mu, \nu\mu^{-1}\delta\mu)} (a\bowtie b)\\
 &=& \Delta_{(\alpha\gamma\alpha^{-1}, \alpha\beta^{-1}\delta\gamma^{-1}\beta\gamma\alpha^{-1}),
    (\alpha\mu\alpha^{-1}, \alpha\beta^{-1}\nu\mu^{-1}\beta\mu\alpha^{-1})}
     \big(a\circ\beta\alpha^{-1}\bowtie \alpha\mu^{-1}\gamma^{-1}\beta^{-1}\gamma\mu(b) \big) \\
 &=& \Delta^{cop}(a\circ\beta\alpha^{-1})
    \Big(\alpha\mu\alpha^{-1} \otimes \alpha\mu^{-1}\alpha^{-1}\cdot \alpha\beta^{-1}\delta\gamma^{-1}\beta\gamma\alpha^{-1} \cdot \alpha\mu\alpha^{-1} \Big)\Delta(\alpha\mu^{-1}\gamma^{-1}\beta^{-1}\gamma\mu(b)) \\
 &=& \Delta^{cop}(a\circ\beta\alpha^{-1})
    \Big(\alpha\gamma^{-1}\beta^{-1}\gamma\mu \otimes \alpha\mu^{-1}\beta^{-1}\delta\mu \Big)\Delta(b)\\
 &=& (\xi_{(\alpha, \beta)}^{(\gamma, \delta)} \otimes \xi_{(\alpha, \beta)}^{(\mu, \nu)})
     \Big(\Delta^{cop}(a) (\mu \otimes \mu^{-1}\delta\mu) \Delta(b)\Big) \\
 &=& (\xi_{(\alpha, \beta)}^{(\gamma, \delta)} \otimes \xi_{(\alpha, \beta)}^{(\mu, \nu)})\circ
     \Delta_{(\gamma, \delta), (\mu, \nu)}(a\bowtie b).
 \end{eqnarray*}

 It is easy to check that $\varepsilon_{D(A, B)} \circ \xi_{(\alpha, \beta)}^{(\iota, \iota)} = \varepsilon_{D(A, B)}$
 for any $(\alpha, \beta)\in G$.

 Finally, we need to check that $\xi_{(\alpha, \beta)} \circ \xi_{(\gamma, \delta)} = \xi_{(\alpha, \beta)\ast (\gamma, \delta)}$.
 Let $a\bowtie b\in A\bowtie B_{(\mu, \nu)}$, wo do the calculations as follows:
 \begin{eqnarray*}
 && \xi_{(\alpha, \beta)}^{(\gamma, \delta)\ast (\mu, \nu)\ast(\gamma, \delta)^{-1}}
    \Big( \xi_{(\gamma, \delta)}^{(\mu, \nu)} (a\bowtie b) \Big) \\
 &=& \xi_{(\alpha, \beta)}^{(\gamma\mu\gamma^{-1}, \gamma\delta^{-1}\nu\mu^{-1}\delta\mu\gamma^{-1})}
    \Big(a\circ \delta\gamma^{-1} \bowtie \gamma\mu^{-1}\delta\mu(b) \Big) \\
 &=& a\circ \delta\gamma^{-1}\beta^{-1}\alpha \bowtie \alpha\cdot \gamma\mu^{-1}\gamma^{-1} \cdot \beta^{-1} \cdot \gamma\mu\gamma^{-1}
    \gamma\mu^{-1}\delta^{-1}\mu(b) \\
 &=& a\circ \delta\gamma^{-1}\beta^{-1}\alpha \bowtie \alpha \gamma\mu^{-1}\gamma^{-1}\beta^{-1} \gamma\delta^{-1}\mu(b) \\
 &=& \xi_{(\alpha\gamma, \delta\gamma^{-1}\beta\gamma)}^{(\mu, \nu)} (a\bowtie b)
 = \xi_{(\alpha, \beta)\ast (\gamma, \delta)} (a\bowtie b).
 \end{eqnarray*}
 Therefore $\xi: G \longrightarrow Aut(\mathcal{A})$ is a crossing action.
 $\hfill \blacksquare$
 \\

 From Theorem \thesection.2 and Proposition \thesection.3, we get that
 $D(A, B) = \bigoplus_{(\alpha, \beta)\in G} A \bowtie B_{(\alpha, \beta)}$
 is a  multiplier Hopf $T$-coalgebra introduced in \cite{YZM13}.
 \\

 Recall from \cite{DV04} that a canonical multiplier $W$ for $\langle A, B\rangle$ ia an invertible element in $M(B\otimes A)$ such that
 $\langle W, a\otimes b\rangle = \langle a, b\rangle$ for all $a\in A$ and $b\in B$.
 Observe that we use the extension of the non-degenerate bilinear form $\langle B\otimes A, A\otimes B\rangle$ to
 $\langle M(B\otimes A), A\otimes B\rangle$. If there is a canonical multiplier $W$ in $M(B\otimes A)$, then it is unique.
 Similar to Proposition 4.4 in \cite{DV04}, we have $(\Delta_{B} \otimes \iota_{A})W = W^{13} W^{23}$
 and $(\iota_{B} \otimes\Delta_{A})W = W^{12} W^{13}$.
 \\

 \textbf{Lemma \thesection.2} Let $W$ be the canonical multiplier in $M(B\otimes A)$. Then
 \begin{itemize}
 \item[(1)] in $M\big( A \bowtie B_{(\alpha, \beta)} \otimes A \big)$,
 \begin{eqnarray}
 (\beta^{-1} \otimes \iota)(W) \Delta^{cop}(a) = \Big(\Delta(a)\circ (\iota \otimes \alpha\beta^{-1})\Big)(\beta^{-1} \otimes \iota)(W). \label{c}
 \end{eqnarray}

 \item[(2)] in $M\big(B \otimes A \bowtie B_{(\gamma, \delta)} \big)$,
 \begin{eqnarray}
 (\beta^{-1} \otimes \iota)(W) (\gamma\otimes \gamma^{-1}\beta\gamma)\Delta(b)
 = (\beta^{-1}\delta\gamma^{-1}\beta\gamma\otimes \gamma^{-1}\beta\gamma)\Delta^{cop}(b)(\beta^{-1} \otimes \iota)(W). \label{d}
 \end{eqnarray}
 \end{itemize}

 \emph{Proof}
 We prove (1). The proof of (2) is similar. We claim that in the multiplier algebra $M(A \bowtie B_{(\alpha, \beta)})$
 \begin{eqnarray*}
 && \Big(\iota \otimes\langle\cdot, b\rangle \Big)\Big( (\beta^{-1} \otimes \iota)(W) \Delta^{cop}(a) \Big) (x\bowtie y) \\
 &=& \Big(\iota \otimes \langle\cdot, b\rangle \Big)
    \Big[ \Big(\Delta(a)\circ (\iota \otimes \alpha\beta^{-1})\Big)(\beta^{-1} \otimes \iota)(W) \Big] (x\bowtie y)
 \end{eqnarray*}
 for all $b\in B$ and $x\bowtie y\in A \bowtie B_{(\alpha, \beta)}$.

 The left-hand side of the above claim is given by
 \begin{eqnarray*}
 && \Big(\iota \otimes\langle\cdot, b\rangle \Big)\Big( (\beta^{-1} \otimes \iota)(W) \Delta^{cop}(a) \Big) (x\bowtie y) \\
 &=& (\iota \otimes\langle\cdot, b_{(1)}\rangle)\Big((\beta^{-1} \otimes \iota)W \Big)
     \langle a_{(1)}, b_{(2)}\rangle (a_{(2)}x \bowtie y) \\
 &=& \langle a_{(1)}, b_{(2)}\rangle \big(1\bowtie \beta^{-1}(b_{(1)}) \big)(a_{(2)}x \bowtie y).
 \end{eqnarray*}
 Take $a'\in A$ such that $b=b\blacktriangleleft a'$. Then The right-hand side of the above claim is given by
 \begin{eqnarray*}
 && \Big(\iota \otimes\langle\cdot, b\rangle \Big)
    \Big[ \Big(\Delta(a)\circ (\iota \otimes \alpha\beta^{-1})\Big)(\beta^{-1} \otimes \iota)(W) \Big] (x\bowtie y) \\
 &=& \Big(\iota \otimes\langle\cdot, b\rangle \Big)
     \Big[ \Big(a_{(1)}\bowtie 1_{M(B)} \otimes a' \big( a_{(2)}\circ \alpha\beta^{-1}\big) \Big)(\beta^{-1} \otimes \iota)(W) \Big] (x\bowtie y) \\
 &=& \langle a' \big( a_{(2)}\circ \alpha\beta^{-1}\big), b_{(1)}\rangle (a_{(1)}\bowtie 1_{M(B)})
    (\langle\cdot, b_{(2)}\rangle \otimes \iota) \Big((\beta^{-1} \otimes \iota)W \Big) (x \bowtie y) \\
 &=& \langle  a_{(2)}\circ \alpha\beta^{-1}, b_{(1)}\rangle \big( a_{(1)}\bowtie \beta^{-1}(b_{(2)}) \big)  (x \bowtie y) \\
 &=& \langle  a_{(2)}, \alpha\beta^{-1}(b_{(1)})\rangle \big( a_{(1)}\bowtie \beta^{-1}(b_{(2)}) \big)  (x \bowtie y).
 \end{eqnarray*}
 Following the commutation rule (\ref{b}) we obtain that the claim is proven.
 Now we get the assertion (1) by using the facts that the pairing is a non-degenerate bilinear form
 and that the product in $A \bowtie B_{(\alpha, \beta)}$ is non-degenerate.
 $\hfill \blacksquare$
 \\

 \textbf{Theorem \thesection.3}
 Let $A$ and $B$ be regular multiplier Hopf algebras,
 $\langle A, B\rangle$ be the multiplier Hopf algebras pairing with the canonical multiplier $W$.
 Then $D(A, B) = \bigoplus_{(\alpha, \beta)\in G} A \bowtie B(\alpha, \beta)$
 is quasitriangular  with a generalized R-matrix given by
 \begin{eqnarray*}
 R = \sum_{(\alpha, \beta), (\gamma, \delta)\in G} R_{(\alpha, \beta), (\gamma, \delta)}
   = \sum_{(\alpha, \beta), (\gamma, \delta)\in G} (\beta^{-1} \otimes \iota)(W).
 \end{eqnarray*}

 \emph{Proof}
 By Proposition 3.2 $(\beta^{-1} \otimes \iota)(W)$ can be embedded in $M(A \bowtie B(\alpha, \beta) \otimes A \bowtie B(\gamma, \delta))$
 by $b\otimes a \hookrightarrow 1_{M(A)}\bowtie b\otimes a \bowtie 1_{M(B)}$.
 Hence, $R_{(\alpha, \beta), (\gamma, \delta)}$ is an element in $M(A \bowtie B_{(\alpha, \beta)}\otimes A \bowtie B_{(\gamma, \delta)})$.
 In the following, we need to check four axioms of quasitriangular structure.

 Firstly, it is easy to check that $\big(\xi_{(\mu, \nu)} \otimes \xi_{(\mu, \nu)} \big)R = R$, since
 \begin{eqnarray*}
 \big(\xi_{(\mu, \nu)} \otimes \xi_{(\mu, \nu)} \big) R_{(\alpha, \beta), (\gamma, \delta)}
 &=&  \big(\xi_{(\mu, \nu)} \otimes \xi_{(\mu, \nu)} \big) (\beta^{-1} \otimes \iota)(W) \\
 &=& \Big(\mu\alpha^{-1}\nu^{-1}\alpha\beta^{-1} \otimes (\cdot)\circ \nu\mu^{-1}\Big) (W), \\
 R_{(\mu, \nu)\ast(\alpha, \beta)\ast(\mu, \nu)^{-1}, (\mu, \nu)\ast (\gamma, \delta) \ast (\mu, \nu)^{-1}}
 &=& R_{(\mu\alpha\mu^{-1}, \mu\nu^{-1}\beta\alpha^{-1}\nu\alpha\mu^{-1}), (\mu\gamma\mu^{-1}, \mu\nu^{-1}\delta\gamma^{-1}\nu\gamma\mu^{-1})} \\
 &=& \big(\mu\alpha^{-1}\nu^{-1}\alpha\beta^{-1}\nu\mu^{-1} \otimes \iota\big) (W)
 \end{eqnarray*}
 and $\big(\iota\otimes (\cdot)\circ \alpha\big)W = (\alpha\otimes\iota)W$.

 Secondly, we need to check that
 \begin{eqnarray*}
 &&(\Delta_{(\alpha, \beta), (\gamma, \delta)} \otimes \iota) R_{(\alpha, \beta)\ast(\gamma, \delta), (\mu, \nu)}
 = \Big((\iota\otimes \xi_{(\gamma, \delta)^{-1}})R_{(\alpha, \beta), (\gamma, \delta) \ast(\mu, \nu) \ast (\gamma, \delta)^{-1}}\Big)_{13}
     \Big( R_{(\gamma, \delta), (\mu, \nu)} \Big)_{23},  \\
 &&(\iota \otimes\Delta_{(\gamma, \delta), (\mu, \nu)}) R_{(\alpha, \beta), (\gamma, \delta)\ast(\mu, \nu)}
 = \Big( R_{(\alpha, \beta), (\mu, \nu)} \Big)_{13}
     \Big( R_{(\alpha, \beta), (\gamma, \delta)} \Big)_{12}.
 \end{eqnarray*}
 We only check the first equation, the second one is similar.
 \begin{eqnarray*}
 && (\Delta_{(\alpha, \beta), (\gamma, \delta)} \otimes \iota) R_{(\alpha, \beta)\ast(\gamma, \delta), (\mu, \nu)}
 = (\Delta_{(\alpha, \beta), (\gamma, \delta)} \otimes \iota) (\beta^{-1} \otimes \iota)(W) \\
 &=& (\beta^{-1}\gamma\delta^{-1}\otimes \delta^{-1}\otimes \iota)(\Delta_{B}\otimes \iota)(W)
 = (\beta^{-1}\gamma\delta^{-1}\otimes \delta^{-1}\otimes \iota)(W^{13}W^{23}),
 \end{eqnarray*}
 and
 \begin{eqnarray*}
 && \Big((\iota\otimes \xi_{(\gamma, \delta)^{-1}})R_{(\alpha, \beta), (\gamma, \delta) \ast(\mu, \nu) \ast (\gamma, \delta)^{-1}}\Big)_{13}
     \Big( R_{(\gamma, \delta), (\mu, \nu)} \Big)_{23} \\
 &=& \Big((\iota\otimes \xi_{(\gamma, \delta)^{-1}})R_{(\alpha, \beta),
     (\gamma\mu\gamma^{-1}, \gamma\delta^{-1}\nu\mu^{-1}\delta\mu\gamma^{-1})}\Big)_{13} \Big( R_{(\gamma, \delta), (\mu, \nu)} \Big)_{23} \\
 &=& \Big((\iota\otimes \xi_{(\gamma^{-1}, \gamma\delta^{-1}\gamma^{-1})})((\beta^{-1} \otimes \iota)(W))\Big)_{13}
     \Big( (\delta^{-1} \otimes \iota)(W) \Big)_{23} \\
 &=& \Big((\beta^{-1} \otimes (\cdot)\circ \gamma\delta^{-1})(W) \Big)_{13} \Big( (\delta^{-1} \otimes \iota)(W) \Big)_{23} \\
 &=& \Big((\beta^{-1}\gamma\delta^{-1} \otimes \iota)(W) \Big)_{13} \Big( (\delta^{-1} \otimes \iota)(W) \Big)_{23} \\
 &=& (\beta^{-1}\gamma\delta^{-1}\otimes \delta^{-1}\otimes \iota)(W^{13}W^{23}).
 \end{eqnarray*}

 Finally, we will check the last axiom:
 \begin{eqnarray*}
 R_{(\alpha, \beta), (\gamma, \delta)}\Delta_{(\alpha, \beta), (\gamma, \delta)}(a\bowtie b) = (\widetilde{\Delta}_{(\alpha, \beta)\ast(\gamma, \delta)\ast(\alpha, \beta)^{-1}, (\alpha, \beta)})^{cop}(a\bowtie b)R_{(\alpha, \beta), (\gamma, \delta)}.
 \end{eqnarray*}
 By Lemma \thesection.2, on the one hand,
 \begin{eqnarray*}
 &&  R_{(\alpha, \beta), (\gamma, \delta)} \Delta_{(\alpha, \beta), (\gamma, \delta)}(a\bowtie b) \\
 &=& (\beta^{-1} \otimes \iota)(W) \Delta^{cop}(a)(\gamma\otimes \gamma^{-1}\beta\gamma)\Delta(b) \\
 &\overset{(\ref{c})}{=}&
     \Big(\Delta(a)\circ (\iota \otimes \alpha\beta^{-1})\Big)(\beta^{-1} \otimes \iota)(W)
     (\gamma\otimes \gamma^{-1}\beta\gamma)\Delta(b) \\
 &\overset{(\ref{d})}{=}&
     \Big(\Delta(a)\circ (\iota \otimes \alpha\beta^{-1})\Big)
     (\beta^{-1}\delta\gamma^{-1}\beta\gamma\otimes \gamma^{-1}\beta\gamma)\Delta^{cop}(b)(\beta^{-1} \otimes \iota)(W),
 \end{eqnarray*}
 on the other hand,
 \begin{eqnarray*}
 && (\widetilde{\Delta}_{(\alpha, \beta)\ast(\gamma, \delta)\ast(\alpha, \beta)^{-1}, (\alpha, \beta)})^{cop}(a\bowtie b)
     R_{(\alpha, \beta), (\gamma, \delta)} \\
 &=& \Big[ \tau
     \Big( \xi_{(\alpha, \beta)^{-1}}^{(\alpha, \beta)\ast(\gamma, \delta)\ast(\alpha, \beta)^{-1}} \otimes \iota\Big)
     \Delta_{(\alpha, \beta)\ast(\gamma, \delta)\ast(\alpha, \beta)^{-1}, (\alpha, \beta)}  (a\bowtie b) \Big]
     R_{(\alpha, \beta), (\gamma, \delta)} \\
 &=& \Big[ \tau
     \Big( \xi_{(\alpha, \beta)^{-1}}^{(\alpha\gamma\alpha^{-1}, \alpha\beta^{-1}\delta\gamma^{-1}\beta\gamma\alpha^{-1})} \otimes \iota\Big)
     \Delta_{(\alpha\gamma\alpha^{-1}, \alpha\beta^{-1}\delta\gamma^{-1}\beta\gamma\alpha^{-1}), (\alpha, \beta)}  (a\bowtie b) \Big]
     (\beta^{-1} \otimes \iota)(W) \\
 &=& \Big[ \tau
     \Big( \xi_{(\alpha, \beta)^{-1}}^{(\alpha\gamma\alpha^{-1}, \alpha\beta^{-1}\delta\gamma^{-1}\beta\gamma\alpha^{-1})} \otimes \iota\Big)
     \Delta^{cop}(a)(\alpha \otimes \beta^{-1}\delta\gamma^{-1}\beta\gamma)\Delta(b) \Big]
     (\beta^{-1} \otimes \iota)(W) \\
 &=& \Big[ \Big(\iota \otimes \xi_{(\alpha, \beta)^{-1}}^{(\alpha\gamma\alpha^{-1}, \alpha\beta^{-1}\delta\gamma^{-1}\beta\gamma\alpha^{-1})}\Big)
     \Delta(a)(\beta^{-1}\delta\gamma^{-1}\beta\gamma \otimes\alpha)\Delta^{cop}(b) \Big] (\beta^{-1} \otimes \iota)(W) \\
 &=& \Big(\Delta(a)\circ (\iota \otimes \alpha\beta^{-1})\Big)
     (\beta^{-1}\delta\gamma^{-1}\beta\gamma\otimes \gamma^{-1}\beta\gamma)\Delta^{cop}(b)(\beta^{-1} \otimes \iota)(W).
 \end{eqnarray*}
 Thus $R$ is a quasitriangular structure in $D(A, B)$.
 $\hfill \blacksquare$
 \\

 \textbf{Example \thesection.5}
 With the notations as Example 3.4. Then by Theorem \thesection.3 the quasitriangular structure is given by
 \begin{eqnarray*}
 R
 &=& \sum_{(\alpha, \beta), (\gamma, \delta)\in G} R_{(\alpha, \beta), (\gamma, \delta)} \\
 &=& \sum_{(\alpha, \beta), (\gamma, \delta)\in G} \sum_{g, h\in H}
      \Big(1_{\widehat{D(H)}} \bowtie (\delta_{\beta^{-1}g\beta} \propto \beta^{-1}h\beta)\Big)\otimes
      \Big((g \propto \delta_{h}) \bowtie 1_{D(H)}\Big).
 \end{eqnarray*}

\section{Applications to Hopf algebras}
\def\theequation{\thesection.\arabic{equation}}
\setcounter{equation}{0}

 In this section, we apply our results as above to the usual Hopf algebras and derive some interesting results.
 First let $H$ be a coFrobenius Hopf algebra with a left integral $\varphi$,
 then by \cite{Zh99} $\widehat{H} = \varphi(\cdot H)$ is a regular multiplier Hopf algebra with integrals,
 and $\langle\widehat{H}, H\rangle$ is a multiplier Hopf pairing. 
 Then by Theorem 3.3 we obtain the following result, which give a positive answer to the question in the introduction.
 \\

 \textbf{Theorem \thesection.1}
 Let $H$ be a coFrobenius Hopf algebra with its dual multiplier Hopf algebra $\widehat{H}$.
 Then $D(\widehat{H}, H) = \bigoplus_{(\alpha, \beta)\in G} \widehat{H} \bowtie H_{(\alpha, \beta)}$
 is a $G$-cograded multiplier Hopf algebra with the following strucrures:
 \begin{itemize}
 \item For any $(\alpha, \beta)\in G$, $\widehat{H} \bowtie H_{(\alpha, \beta)}$ has the multiplication given by
 \begin{eqnarray*}
 (p\bowtie h)(q\bowtie l) = p \big( \alpha(h_{(1)}) \blacktriangleright q \blacktriangleleft S^{-1}\beta(h_{(3)})  \big) \bowtie h_{(2)}l
 \end{eqnarray*}
 for $p, q\in \widehat{H}$ and $h, l\in H$.

 \item The comultiplication on $D(\widehat{H}, H)$ is given by:
 \begin{eqnarray*}
 && \Delta_{(\alpha, \beta), (\gamma, \delta)}:
 \widehat{H} \bowtie H_{(\alpha, \beta)\ast (\gamma, \delta)}
 \longrightarrow \widehat{H} \bowtie H_{(\alpha, \beta)} \otimes \widehat{H} \bowtie H_{(\gamma, \delta)}, \\
 && \Delta_{(\alpha, \beta), (\gamma, \delta)} (p\bowtie h) = \Delta^{cop}(p)(\gamma\otimes \gamma^{-1}\beta\gamma)\Delta(h).
 \end{eqnarray*}

 \item The counit $\varepsilon_{D(\widehat{H}, H)} = \varepsilon_{\widehat{H}}\otimes \varepsilon_{H}$.

 \item For any $(\alpha, \beta)\in G$, the antipode is given by
 \begin{eqnarray*}
 && S: \widehat{H} \bowtie H_{(\alpha, \beta)} \longrightarrow \widehat{H} \bowtie H_{(\alpha, \beta)^{-1}}, \\
 && S_{(\alpha, \beta)}(p\bowtie h) = T(\alpha\beta S(h) \otimes S^{-1}(p))
    \mbox{ in } \widehat{H} \bowtie H_{(\alpha, \beta)^{-1}}.
 \end{eqnarray*}
 \end{itemize}

 If furthermore there is a cointegral $t\in A$ such that $\varphi(t)=1$.
 Then by Theorem 4.3 $D(\widehat{H}, H) = \bigoplus_{(\alpha, \beta)\in G} \widehat{H} \bowtie H_{(\alpha, \beta)}$
 admits a quasitriangular structure.
 \\

 \textbf{Theorem \thesection.2}
 Let $H$ be a coFrobenius Hopf algebra with its dual multiplier Hopf algebra $\widehat{H}$.
 Then $\mathcal{A} = \bigoplus_{(\alpha, \beta)\in G} \widehat{H} \bowtie H(\alpha, \beta)$
 is a quasitriangular $G$-cograded multiplier Hopf algebra with a generalized R-matrix given by
 \begin{eqnarray*}
 R = \sum_{(\alpha, \beta), (\gamma, \delta)\in G} R_{(\alpha, \beta), (\gamma, \delta)}
   = \sum_{(\alpha, \beta), (\gamma, \delta)\in G} \varepsilon\bowtie\beta^{-1}(t(\cdot \varphi_{(2)})) \otimes S^{-1}(\varphi_{(1)})\bowtie 1.
 \end{eqnarray*}
 \\

 \textbf{Example \thesection.3} Let $H$ be an infinite group with unit $e$.
 We denote by $KH$ the corresponding group algebra and by $K(H)$ the classical dual multiplier Hopf algebra.
 $G = Aut_{Hopf}(H)\times Aut_{Hopf}(H)$ is a group with product (3.1).
 Let $\alpha\in H$, we define $\alpha(h)=\alpha h \alpha^{-1}$.
 Then $\alpha\in Aut_{Hopf}(H)$, and by Theorem \thesection.3 we can construct a $G$-cograded multiplier Hopf algebra $D(K(H), KH)$ with the
 multiplication in $K(H)\bowtie KH_{(\alpha, \beta)}$, comultiplication, counit in $K(H)\bowtie KH_{(\iota, \iota)}$, antipode as follows:
 \begin{eqnarray*}
 && (\delta_{p}\bowtie g)(\delta_{q}\otimes h) = \delta_{p}\delta_{\beta g \beta^{-1} q\alpha g^{-1} \alpha^{-1} }\bowtie gh,\\
 && \Delta_{(\alpha, \beta), (\gamma, \delta)} (\delta_{p}\bowtie h)
   = \sum_{s\in H} \delta_{s^{-1}p}\bowtie \gamma h \gamma^{-1} \otimes \delta_{s}\bowtie \gamma^{-1}\beta\gamma h \gamma^{-1}\beta^{-1}\gamma, \\
 && \varepsilon(\delta_{p}\bowtie g) = \delta_{p, e}, \\
 && S_{(\alpha, \beta)}(\delta_{p}\bowtie h) = \delta_{\alpha h^{-1}\alpha^{-1} p^{-1}\beta h \beta^{-1}}
    \otimes \alpha\beta h^{-1}\beta^{-1}\alpha^{-1}.
 \end{eqnarray*}
 The quasitriangular structure is given by
 \begin{eqnarray*}
 R = \sum_{(\alpha, \beta), (\gamma, \delta)\in G} R_{(\alpha, \beta), (\gamma, \delta)}
   = \sum_{(\alpha, \beta), (\gamma, \delta)\in G; g\in H} \varepsilon\bowtie\beta^{-1} g \beta \otimes \delta_{g}\bowtie e.
 \end{eqnarray*}
 \\

 Let $B=H$ be a finite dimensional Hopf algebra and $A=H^*$ be the dual Hopf algebra. Then we can get the following result,
 which is constructed by Panaite and Staic Mihai in \cite{PS07}.
 \\

  \textbf{Theorem \thesection.4}
 Let $H$ be a finite dimensional Hopf algebra.
 Then $D(H) = \bigoplus_{(\alpha, \beta)\in G} H^* \bowtie H_{(\alpha, \beta)}$
 is a $G$-cograded multiplier Hopf algebra with the following strucrures:
 \begin{itemize}
 \item For any $(\alpha, \beta)\in G$, $H^* \bowtie H_{(\alpha, \beta)}$ has the multiplication given by
 \begin{eqnarray*}
 (p\bowtie h)(q\bowtie l) = p \big( \alpha(h_{(1)}) \blacktriangleright q \blacktriangleleft S^{-1}\beta(h_{(3)})  \big) \bowtie h_{(2)}l
 \end{eqnarray*}
 for $p, q\in H^*$ and $h, l\in H$.

 \item The comultiplication on $D(H^*, H)$ is given by:
 \begin{eqnarray*}
 && \Delta_{(\alpha, \beta), (\gamma, \delta)}:
 H^* \bowtie H_{(\alpha, \beta)\ast (\gamma, \delta)}
 \longrightarrow H^* \bowtie H_{(\alpha, \beta)} \otimes H^* \bowtie H_{(\gamma, \delta)}, \\
 && \Delta_{(\alpha, \beta), (\gamma, \delta)} (p\bowtie h) = \Delta^{cop}(p)(\gamma\otimes \gamma^{-1}\beta\gamma)\Delta(h).
 \end{eqnarray*}

 \item The counit $\varepsilon_{D(H^*, H)} = \varepsilon_{H^*}\otimes \varepsilon_{H}$.

 \item For any $(\alpha, \beta)\in G$, the antipode is given by
 \begin{eqnarray*}
 && S: H^* \bowtie H_{(\alpha, \beta)} \longrightarrow H^* \bowtie H_{(\alpha, \beta)^{-1}}, \\
 && S_{(\alpha, \beta)}(p\bowtie h) = T(\alpha\beta S(h) \otimes S^{-1}(p))
    \mbox{ in } H^* \bowtie H_{(\alpha, \beta)^{-1}}.
 \end{eqnarray*}

 \item The generalized R-matrix given by
 \begin{eqnarray*}
 R = \sum_{(\alpha, \beta), (\gamma, \delta)\in G} R_{(\alpha, \beta), (\gamma, \delta)}
   = \sum_{(\alpha, \beta), (\gamma, \delta)\in G} \varepsilon\bowtie\beta^{-1}(e_{i}) \otimes S^{-1}(e^{i})\bowtie 1,
 \end{eqnarray*}
 where $e_{i}$ and $e^{i}$ are dual basis of $H$ and $H^*$.
 \end{itemize}

\section*{Acknowledgements}

 The work was partially supported by the National Natural Science Foundation of China (Grant No. 11601231),
 the Fundamental Research Fund for the Central Universities (Grant No.KJQN201716) and the Natural Science Foundation of
 Jiangsu Province (Grant No.BK20160708).


\addcontentsline{toc}{section}{References}
\vskip 0.6cm
\begin{thebibliography}{99}



 \bibitem{ADV07} Abd El-Hafez A.T., Delvaux L. and Van Daele, A. (2007). Group-cograded Multiplier Hopf (*-)algebras.
 \emph{Algebras and Representation Theory} 10: 77-95.


 \bibitem{D03} Delvaux, L. (2003). Twisted tensor product of multiplier Hopf (*-) algebras.
 \emph{Journal of Algebra} 269: 285-316.

 \bibitem{DV04}  Delvaux, L. and Van Daele, A. (2004). The Drinfel'd double versus the Heisenberg double for an algebraic quantum group.
 \emph{Journal of Pure and Applied Algebra} 190: 59-84.

 \bibitem{DV07} Delvaux, L. and Van Daele, A. (2007).  The Drinfel¡¯d Double for group-cograded multiplier Hopf algebras.
 \emph{Algebras and Representation Theory} 10(3): 197-221

 \bibitem{DVW05} Delvaux L., Van Daele A. and Wang S. H. (2005). Quasitriangular (G-cograded) multiplier Hopf algebras.
 \emph{Journal of Algebra} 289: 484-514.


 \bibitem{DrV01} Drabant, B. and Van Daele, A. (2001). Pairing and quantum double of multiplier Hopf algebras.
 \emph{Algebras and Representation Theory} 4: 109-132.
%
%
%

 \bibitem{PS07} Panaite F. and Staic Mihai D. (2007). Generalized (anti) Yetter-Drinfel'd modules as components of a braided T-category.
 \emph{Israel Journal of Mathematics} 158: 349-365.

 \bibitem{T} Turaev, V. G. (2000). Homotopy field theory in dimension $3$ and crossed group-categories.
      Preprint GT/0005291

 \bibitem{V94} Van Daele, A. (1994). Multiplier Hopf algebras.
 \emph{Transaction of the American Mathematical Society} 342(2): 917-932.

 \bibitem{V98} Van Daele, A. (1998). An algebraic framework for group duality.
 \emph{Advances in Mathematics} 140(2): 323-366.

 \bibitem{V08} Van Daele, A. (2008). Tools for working with multiplier Hopf algebras.
 \emph{The Arabian Journal for Science and Engineering} 33(2C): 505-527.

%


 \bibitem{YW11a} Yang, T. and Wang, S. H. (2011). A lot of quasitriangular group-cograded multiplier Hopf algebras.
 \emph{Algebras and Representation Theory} 14(5): 959-976.

 \bibitem{YW11} Yang, T. and Wang, S. H. (2011). Constructing new braided $T$-categories over regular multiplier Hopf algebras.
 \emph{Communications in Algebra}, 39(9): 3073-3089.

 \bibitem{YZM13} Yang, T., Zhou, X. and Ma, T. (2013). On braided T-categories over multiplier Hopf algebras.
 \emph{Communications in Algebra}, 41: 2852¨C2868.


 \bibitem{Zh99} Zhang, Y. H. (1999). The quantum double of a coFrobenius Hopf algebra.
 \emph{Communications in Algebra}, 27(3): 1413-1427.

\end{thebibliography}
\end {document}